\title{\vspace{-1cm}A new bound for the cops and robbers problem}
\author{Alex Scott\thanks{Mathematical Institute, 24-29 St Giles', Oxford, OX1 3LB, UK. Email:
{\tt scott@maths.ox.ac.uk}.} \and Benny Sudakov\thanks{Department of Mathematics,
UCLA,  Los Angeles, CA 90095. Email: {\tt bsudakov@math.ucla.edu}.
Research supported in part by NSF CAREER award DMS-0812005 and by
USA-Israeli BSF grant.}}

\documentclass[11pt]{article}
\usepackage{amsfonts,amssymb,amsmath,latexsym}

\newenvironment{proof}
      {\medskip\noindent{\bf Proof.}\hspace{1mm}}
      {\hfill$\Box$\medskip}

\def\qed{\ifvmode\mbox{ }\else\unskip\fi\hskip 1em plus 10fill$\Box$}

\newtheorem{theorem}{Theorem}[section]

\newtheorem{lemma}[theorem]{Lemma}

\newtheorem{claim}[theorem]{Claim}

\begin{document}
\date{}
\maketitle

\begin{abstract}
In this short paper we study the game of cops and robbers, which 
is played on the vertices of some fixed graph $G$. Cops and a robber are allowed to move along the edges of $G$ and the 
goal of cops is to capture the robber.
The cop  number $c(G)$ of $G$ is the minimum number of cops required to win the game.
Meyniel conjectured a long time ago that $O(\sqrt{n})$ cops are enough for any connected $G$ on $n$ vertices.
Improving several previous results, we prove that the cop number of $n$-vertex graph is at most $n 2^{-(1+o(1))\sqrt{\log n}}$.
\end{abstract}

\section{Introduction}
Let $G$ be a simple, undirected, connected graph on $n$ vertices. The game of 
{\em Cops and Robbers}, which was introduced almost thirty years by Nowakowski and Winkler \cite{NW} and by Quilliot \cite{Q}, 
is played on the vertices of $G$ as follows. There are two players, a set of $k \geq 1$ cops and one robber. The game begins by the
cops occupying some set of $k$ vertices of $G$ and then the robber also chooses a vertex. Afterward they move alternatively, first cops and then robber, along the edges of graph $G$. At every step each cop or robber is allowed to move to any neighboring vertex or do nothing and stay where they are. 
Multiple cops are allowed to occupy the same vertex. The cops win if at some time there is a cop at the same vertex as the robber; otherwise, the robber wins. 
The minimum number of cops for which there is a
winning strategy, no matter how robber plays, is called {\em the cop number of $G$} and is denoted by $c(G)$. Note that this number is clearly at most $n$ and also
that the initial position of cops does not matter, since $G$ is connected.

The cop number was introduced by Aigner and Fromme \cite{AF} who proved for example that if $G$ is planar, then $c(G) \leq 3$.
They also observed that if $G$ has has girth at least $5$ (i.e., no cycles of length $\leq 4$) then its cop number is at least the minimum degree of $G$.
In particular, this together with the well known construction of dense graphs of girth $5$ shows that there are $n$-vertex graphs which require at least $\Omega(\sqrt{n})$ cops. A quarter century ago, Meyniel conjectured that this is tight and $O(\sqrt{n})$ cops is always sufficient. The first nontrivial upper 
bound for this problem was obtained by Frankl \cite{F}, who proved that $c(G) \leq O(\frac{n \log \log n}{\log n})$. This was later
improved by Chiniforooshan \cite{Ch} to $O(\frac{n}{\log n})$. All logarithms in this paper are binary.

The game of cops and robbers was also studied by
Andreae \cite{A}, Berardicci and Intrigila \cite{BI}, Alspach \cite{Al} and for random graph by Bollob\'as, Kun and Leader \cite{BKL} and by 
\L uczak and Pra\l at \cite{LP}. The aim of this short paper is to prove the following new bound for the cops and robber problem.

\begin{theorem}
\label{main}
The cop number of any connected $n$-vertex graph is at most $n 2^{-(1+o(1))\sqrt{\log n}}$.
\end{theorem}

\section{Proof of the main result}
The proof of the main theorem has two ingredients. The first is a result of Aigner and Fromme 
which says that one cop can control the shortest path between two vertices. Let $P$ be the shortest path 
in $G$ between two vertices $u$ and $v$. Then the following lemma was proved in \cite{AF}.

\begin{lemma}
\label{shortest}
One cop can move along the vertices of $P$ such that, after finite number of steps, if the robber 
ever visits $P$ then he will be caught in the next step.
\end{lemma}

This result can be used for graphs with large diameter, since we can delete a long path from such a graph and use induction.
The case when $G$ has relatively small diameter will be treated using the following key lemma, which we think has 
independent interest.

\begin{lemma}
\label{expand}
Let $G$ be a graph on $n \geq 2^{30}$ vertices with diameter $D \leq 2^{\sqrt{\log n}} /\log^3 n$. Then 
$c(G) \leq n (\log n)^3 2^{-\sqrt{\log n}}$.
\end{lemma}

\begin{proof}
Let $t=\sqrt{\log n} -3 \log \log n$. Consider random subsets ${\cal C}_1, \ldots, {\cal C}_{t+1}$ of $G$. For every vertex of $G$ we put it into each
${\cal C}_j$ with probability $(\log n)^2 2^{-\sqrt{\log n}}$ (so a single vertex may be in many ${\cal C}_i$). 
Since $|{\cal C}_j|$ is binomially distributed with expectation
$\mu=n (\log n)^2 2^{-\sqrt{\log n}}$, 
by the standard Chernoff-type estimates (see, e.g., Appendix A in \cite{AS}), we have that 
the probability that ${\cal C}_j$ has more than $2\mu=2n (\log n)^2 2^{-\sqrt{\log n}}$ vertices is at most
$e^{-\mu/3} <n^{-2}$. Then, with probability at least $1-\frac{\log n}{n^2}>0.9$, the sum of the sizes of these sets is at most
$n (\log n)^3 2^{-\sqrt{\log n}}$ which will give our final bound. 

For every subset $A$ of $G$ and integer $i$ let $B(A,i)$ be the ball of radius $i$ around $A$, i.e., all the vertices of $G$
which can be reached from some vertex in $A$ by a path of length at most $i$. We need the following simple claim.


\begin{claim}
\label{random}
The following statement holds with probability $0.9$:
for every $A\subset V(G)$ such that $|A| \leq  n 2^{-\sqrt{\log n}}$, every $i \leq t$ such that 
$|B(A,2^i)| \geq 2^{\sqrt{\log n}}|A|$, and every $j$,
$$|B(A,2^i)\cap {\cal C}_j|\ge |A|.$$
\end{claim}

\begin{proof}
Let $|A|=a$. Note that for any fixed $A, i, j$ the number of points from ${\cal C}_j$ in $B(A,2^i)$ is binomial distributed
with expectation at least $a \log ^2 n$. Thus by Chernoff type estimates the probability that it is smaller than $a$ is at most 
$e^{-a\log^2 n/3}$.
The number of sets of size $a$ is $\binom na$ and the number of pairs of indices $i,j$ is at most $\log n$, and the result follows 
since $\log n \sum_a \binom na e^{-a\log^2 n/3} <0.1$.
\end{proof}

Thus we can choose subsets ${\cal C}_1, \ldots, {\cal C}_{t+1}$ that satisfy the assertion of this claim and such that the sum of their sizes
is at most $n (\log n)^3 2^{-\sqrt{\log n}}$. We place at each vertex $u$ of $G$ one cop for each set ${\cal C}_j$ that contains $u$.
We will show that these cops can always catch the robber. 

Suppose that the robber is located at vertex $v$ of the graph $G$. Since the cops move first, note that the degree of this vertex is at most
$2^{\sqrt{\log n}}$. Otherwise by the claim we already have a cop in its neighborhood who will
catch the robber in the first move. Consider the following sequence of sets $A_i, D_i$ for $i=1, \ldots t+1$ defined
recursively. 

Let $A_1$ be a the largest subset of $B(v,1)$ such that $|B(A_1,1)| <2^{\sqrt{\log n}}|A_1|$ and let $D_1=B(v,1) - A_1$.
Note that for any $X \subset D_1$ we have that $|B(X,1)| \geq 2^{\sqrt{\log n}} |X|$ since otherwise we could increase $A_1$ by
adding $X$ to it. Thus by Claim \ref{random} we have that $B(X,1)$ contains more than $|X|$ points from ${\cal C}_1$. Consider an auxiliary bipartite graph with
vertex classes
$D_1$ and ${\cal C}_1$, where we join a vertex $u$ in $D_1$ to $w$ in ${\cal C}_1$ if $w$ is within distance at most one from $u$. Note that, 
by the above discussion, it follows from Hall's Theorem that this bipartite graph has a complete matching from $D_1$ to ${\cal C}_1$. 
Therefore for every vertex $u \in D_1$ we can choose a corresponding vertex $w$ in ${\cal C}_1$, which is within distance one from $u$, and such that all
$w$ are distinct. We can therefore ask the cops in ${\cal C}_1$ to occupy $D_1$ in the first round.

Now suppose we have already defined $A_i, D_i$. Consider the ball $B(A_i, 2^{i-1})$. Let $A_{i+1}$ be the 
largest subset of $B(A_i,2^{i-1})$ such that $|B(A_{i+1},2^i)| <2^{\sqrt{\log n}}|A_{i+1}|$ and let $D_{i+1}=B(A_i,2^{i-1}) - A_{i+1}$.
Again for any subset $X \subset D_{i+1}$ we have that $|B(X,2^i)| \geq 2^{\sqrt{\log n}} |X|$. 
Thus by Lemma \ref{random} we have that $B(X,2^i)$ contains more than $|X|$ points from ${\cal C}_{i+1}$. 
As above, we consider an auxiliary bipartite graph with parts
$D_{i+1}$ and ${\cal C}_{i+1}$ such that a vertex $u$ in $D_{i+1}$ is joined to $w$ in ${\cal C}_{i+1}$ if $w$ is within distance at most $2^i$ from $u$. 
The bipartite graph has a complete matching from $D_{i+1}$ to ${\cal C}_{i+1}$ 
and so for every vertex $u \in D_{i+1}$ we can choose a corresponding vertex $w$ in 
${\cal C}_{i+1}$, which is within distance at most $2^i$ from $u$, and such that all
$w$ are distinct. We will ask the cops in ${\cal C}_{i+1}$ to occupy $D_{i+1}$ by the end of the first $2^i$ rounds (by rounds we mean here that cops will move first, then robber moves and so
on until both cops and robber have each made $2^i$ moves).

Next we claim that there is an index $s \leq t+1$ such that $A_s$ is empty.
Indeed note that by construction we have $|A_{i+1}| \leq |B(A_i, 2^{i-1})| \leq 2^{\sqrt{\log n}}|A_i|$.
Therefore if $A_{t+1}$ is not empty we have that
$$|A_{t+1}| \leq \left(2^{\sqrt{\log n}}\right)^{t+1} \leq n 2^{-2\sqrt{\log n}}.$$
Moreover we also have that $B(A_{t+1},2^t)$ has size at most $2^{\sqrt{\log n}}|A_{t+1}|\ll n$ and this contradicts the assumption that the diameter of $G$
is at most $2^{\sqrt{\log n}}/\log^3 n$.

Let $s \leq t+1$ be such that $A_s$ is empty. We have shown that we can move cops from the sets ${\cal C}_i, 1 \leq i \leq s$, such that after 
$2^{i-1}$ rounds they will occupy the set $D_i$. Also we have that $B(A_{s-1}, 2^{s-2})=D_s$.
Now we claim that after $2^{s-1}$ rounds we have already caught the robber. 

After one round the robber has made one step so he is in $B(v,1)$. Note that $D_1$ is already occupied by cops after the first round so
in order to survive the robber must be in $A_1=B(v,1)-D_1$. Suppose by induction that after $2^{i-1}$ rounds the robber is in the set 
$A_i$. Consider the next $2^{i-1}$ rounds (so in total $2 \cdot 2^{i-1}=2^i$ rounds). After this many additional rounds, the robber is in
$B(A_i, 2^{i-1})$. Since we know that $D_{i+1}$ is occupied by the cops by time $2^i$, the robber must be
in $ B(A_i, 2^{i-1})-D_{i+1}=A_{i+1}$. 
Thus, arguing by induction, the only place the robber can be after
$2^{i-1}$ rounds without having been caught is in the set $A_i$. Since $A_s$ is empty we are done.
This completes the proof of the lemma.
\end{proof}

\vspace{0.1cm}
\noindent
{\bf Proof of Theorem \ref{main}.}\, 
Let $G$ be a connected graph on $n$ vertices. We prove by induction that $c(G) \leq f(n)= 2 n (\log n)^3 2^{-\sqrt{\log n}}$.  Since $c(G) \leq n$, the result holds trivially when $ 2(\log n)^3 2^{-\sqrt{\log n}} \geq 1$. Thus
we can assume that $\log n \geq 400$.  
If the diameter of $G$ is at most $2^{\sqrt{\log n}}/\log^3 n$ then we are immediately done by Lemma \ref{expand}. Otherwise $G$ contains two vertices such that the shortest path $P$ between them has length at least $D=2^{\sqrt{\log n}}/\log^3 n$. Put one cop on this path. By Lemma \ref{shortest}, after finite number of steps this cop can prevent the robber entering path $P$. Thus we can continue the game on the graph obtained from $G$ by deleting $P$. If this graph is disconnected we continue playing on the connected component which contains the robber. Thus $c(G) \leq 1 +c(G-P)$. By induction we know that
\begin{equation}
\label{eq1}
c(G-P) \leq f(n-D) \leq 2(n-D)(\log n)^3 2^{-\sqrt{\log (n-D)}} \leq	 2 n (\log n)^3 2^{-\sqrt{\log (n-D)}} - 2.
\end{equation}
Since $\log n \geq 400$,  we have $D \leq 2^{-20} n$. Then $\log (n-D) \geq \log n -2D/n$ and therefore
$$\sqrt{\log (n-D)} \geq \sqrt{\log n -2D/n} \geq \sqrt{\log n}- 2D/(n \sqrt{\log n}).$$ 
We also have that $2^{2D/(n \sqrt{\log n})} \leq \big(1+2D/(n \sqrt{\log n})\big)$. Substituting this into
(\ref{eq1}) we get that 
\begin{eqnarray*} 
c(G-P) &\leq& 2 n (\log n)^3 2^{-\sqrt{\log (n-D)}} - 2 \leq 2 n (\log n)^3 2^{-\sqrt{\log n}+2D/(n \sqrt{\log n})} - 2\\
&\leq& 2 n (\log n)^3 2^{-\sqrt{\log n}} \big(1+2D/(n \sqrt{\log n})\big)-2 = f(n)+\frac{4}{\sqrt{\log n}}-2 <f(n)-1.
\end{eqnarray*}
Therefore $c(G) \leq 1+c(G-P)<1+ f(n)-1=f(n)$ which completes the induction step.
\hfill $\Box$

\vspace{0.1cm}

Finally, let us note that the strategy in Lemma \ref{expand} does not 
require the cops to see the robber's movements after the first round.  
The cops therefore have the following randomized strategy against an invisible robber 
(see \cite{IK} for a discussion of cops-and-robbers games where there is a 
cop with limited vision). The cops move to their assigned 
positions, then guess the position of the robber and follow the 
strategy of Lemma \ref{expand}, repeating this process until the robber is 
caught.  This algorithm catches the robber in polynomial expected time.
On a general graph, we can add a further cop on each path $P$ that we have deleted in the proof of Theorem 1.1: 
if these cops move at random on their paths, the robber
will be caught in finite (although not necessarily polynomial) expected time.

\vspace{0.2cm}
\noindent
{\bf Note added in proof.}\, When this paper was written we learn that 
a similar result independently and slightly before us was 
also obtained by Lu and Peng \cite{Lu}. We'd like to thank Krivelevich for bringing 
their paper to our attention.

\end{document}